\documentclass[11pt]{article}
\usepackage{mathrsfs}
\usepackage{epic,eepic,epsf,epsfig}
\usepackage{amsfonts,srcltx,mathrsfs}
\usepackage{amsmath}
\usepackage{amssymb}
\usepackage{amsbsy}
\usepackage{graphicx}
\usepackage{amsfonts}%,srcltx}
\usepackage{color}
\usepackage{mathrsfs}
\usepackage{epic,eepic,epsf,epsfig}
\usepackage{graphicx}
\usepackage{amsfonts,srcltx,mathrsfs}
\usepackage{array}
\usepackage{amsthm}
\newtheorem{theorem}{Theorem}
\newtheorem{lem}{Lemma}
\newtheorem{cor}{Corollary}

\def\f{\noindent}

\begin{document}

\markboth{ et. al}{Two kinds of generalized connectivity of dual cubes}

\title{Two kinds of generalized connectivity of dual cubes}

\author{Shu-Li Zhao$^{1}$, Rong-Xia Hao$^{1}$\footnote{Corresponding author. Email:17118434@bjtu.edu.cn(S.-L. Zhao), rxhao@bjtu.edu.cn (R.-X. Hao), echeng@oakland.edu(E. Cheng)}, Eddie Cheng$^{2}$\\[0.2cm]
{\em\small $^{1}$Department of Mathematics, Beijing Jiaotong University,}\\ {\small\em
Beijing 100044, P.R. China}\\ {\em\small $^{2}$Department of Mathematics and Statistics, Oakland University,}\\ {\small\em Rochester, MI 48309, USA}}

\date{}
\maketitle

Let $S\subseteq V(G)$ and $\kappa_{G}(S)$ denote the maximum number $k$ of edge-disjoint trees $T_{1}, T_{2}, \cdots, \\T_{k}$ in $G$ such that $V(T_{i})\bigcap V(T_{j})=S$ for any $i, j \in \{1, 2, \cdots, k\}$ and $i\neq j$. For an integer $r$ with $2\leq r\leq n$, the {\em generalized $r$-connectivity} of a graph $G$ is defined as $\kappa_{r}(G)= min\{\kappa_{G}(S)|S\subseteq V(G)$ and $|S|=r\}$. The $r$-component connectivity $c\kappa_{r}(G)$ of a non-complete graph $G$ is the minimum number of vertices whose deletion results in a graph with at least $r$ components. These two parameters are both generalizations of traditional connectivity. Except hypercubes and complete bipartite graphs, almost all known $\kappa_{r}(G)$ are about $r=3$. In this paper, we focus on $\kappa_{4}(D_{n})$ of dual cube $D_{n}$. We first show that $\kappa_{4}(D_{n})=n-1$ for $n\geq 4$. As a corollary, we obtain $\kappa_{3}(D_{n})=n-1$ for $n\geq 4$. Furthermore, we show that $c\kappa_{r+1}(D_{n})=rn-\frac{r(r+1)}{2}+1$ for $n\geq 2$ and $1\leq r \leq n-1$.

\medskip

\f {\em Keywords:} Generalized connectivity; Component connectivity; Fault-tolerance; Dual cube.

\section{Introduction}
For an interconnection network, one major concern is the fault tolerance. An interconnection network is usually modelled by a connected graph $G=(V, E)$, where nodes represent processors and edges represent communication links between processors. The connectivity is one of the important parameters to evaluate the reliability and fault tolerance of a network. {\it The connectivity $\kappa (G)$} of a graph $G$ is defined as the minimum number of vertices whose deletion results in a disconnected graph. Whitney~\cite{w} provides another definition of connectivity. For any subset $S=\{u, v\}\subseteq V(G)$, let $\kappa_{G}(S)$ denote the maximum number of internally disjoint paths between $u$ and $v$ in $G$. Then $\kappa(G)=min\{\kappa_{G}(S)|S\subseteq V(G)$ and $|S|=2\}$. As generalizations of the traditional connectivity, the {\it generalized $k$-connectivity} and the $r$-component connectivity were introduced by Chartrand $et$ $al.$~\cite{c} in $1984$.

Let $G$ be a non-complete graph. A $r$-component cut of $G$ is a set of vertices whose deletion results in a graph with at least $r$ components. The $r$-component connectivity $c\kappa_{r}(G)$ of a graph $G$ is the size of the smallest $r$-component cut of $G$. By the definition of $c\kappa_{r}(G)$, it can be seen that $ c\kappa_{r+1}(G)\geq c\kappa_{r}(G)$ for every positive integer $r$. The $r$-component connectivity is an extension of the usual connectivity $c\kappa_{2}(G)$. In \cite{Hsu}, Hsu $et$ $al.$ determined the $r$-component connectivity of the hypercube $Q_{n}$ for $r=2,3,\cdots, n+1$. In \cite{Z}, Zhao $et$ $al.$ determine the $r$-component connectivity of the hypercube $Q_{n}$ for $r=n+2, n+3, \cdots, 2n-4$.

Let $S\subseteq V(G)$ and $\kappa_{G}(S)$ denote the maximum number $k$ of edge-disjoint trees $T_{1}, T_{2}, \cdots, \\T_{k}$ in $G$ such that $V(T_{i})\bigcap V(T_{j})=S$ for any $i, j \in \{1, 2, \cdots, k\}$ and $i\neq j$. For an integer $r$ with $2\leq r\leq n$, the {\em generalized $r$-connectivity} of a graph $G$ is defined as $\kappa_{r}(G)= min\{\kappa_{G}(S)|S\subseteq V(G)$ and $|S|=r\}$. It is a parameter that can measure the reliability of a network $G$ to connect any $k$ vertices in $G$. If $\kappa_{r}(G)=m$, then there are $m$ internally disjoint trees in $G$ and each of them connecting the vertices of $S$, where $|S|=r$ be any $S$ of $V(G)$. The generalized $2$-connectivity is exactly the traditional connectivity. In~\cite{l4}, Li $et$ $al.$ derived that it is NP-complete for a general  graph $G$ to decide whether there are $k$ internally disjoint trees connecting $S$, where $k$ is a fixed integer and $S\subseteq V(G).$ There are some results~\cite{l2,l5} about the upper and lower bounds of the generalized connectivity. In addition, there are some results of the generalized $k$-connectivity for some classes of graphs and most of them are about $k=3$. For example, Chartrand {\em et al.}~\cite{ch} studied the generalized connectivity of complete graphs; Li $et$ $al.$~\cite{l1} first studied the generalized $3$-connectivity of Cartesian product graphs, then Li $et$ $al.$ ~\cite{l8} studied the generalized $3$-connectivity of graph products; Li $et$ $al.$~\cite{l3} studied the generalized connectivity of the complete bipartite graphs and Lin $et$ $al.$~\cite{L} studied the generalized $4$-connectivity of hypercubes. As the Cayley graph has some attractive properties to design interconnection networks, Li $et$ $al.$~\cite{l6} studied the generalized $3$-connectivity of star graphs and bubble-sort graphs and Li $et$ $al.$~\cite{l7} studied the generalized $3$-connectivity of the Cayley graph generated by trees and cycles. Except hypercubes and complete bipartite graphs, almost all known exact value of the generalized $k$-connectivity are about $k=3$.

In this paper, we first show that $\kappa_{4}(D_{n})=\kappa_{3}(D_{n})=n-1$ for $n\geq 4$, then we show that $c\kappa_{r+1}(D_{n})=rn-\frac{r(r+1)}{2}+1$ for $n\geq 2$ and $1\leq r \leq n-1$.

The paper is organized as follows. In section 2, some notations and definitions are given. In section 3, the generalized $4$-connectivity and the generalized $3$-connectivity of dual cubes are determined. In section 4, the $r$-component connectivity of the dual cube for $n\geq 2$ and $2\leq r \leq n-1$ is determined. In section 5, the paper is concluded.

\section{Preliminary}

Let $G=(V, E)$ be a simple, undirected graph. Let $|V(G)|$ and $|E(G)|$ denote the order and size of graph $G$, respectively. Let $V^{\prime}\subseteq V(G)$, then $G[V^{\prime}]$ is the subgraph of $G$ whose vertex set is $V^{\prime}$ and whose edge set consists of all edges of $G$ which have both ends in $V^{\prime}$. For a vertex $v\in V(G)$, the set of neighbours of a vertex $v$ in a graph $G$ is denoted by $N_{G}(v)$. For a vertex set $U \subseteq V(G)$, the {\it neighborhood} of $U$ in $G$ is defined as $N_G(U)=\bigcup\limits_{v\in U}N_{G}(v)-U$. Let $d_{G}(v)$ denote the number of edges incident with $v$ and $\delta(G)$ denote the {\em minimum degree} of the graph $G$. A graph is said to be {\em $k$-regular} if for any vertex $v$ of $G$, $d_{G}(v)=k$. Two $xy$- paths $P$ and $Q$ in $G$ are {\em internally disjoint} if they have no common internal vertices, that is $V(P)\bigcap V(Q)=\{x, y\}$. Let $Y\subseteq V(G)$ and $X\subset V(G)\setminus Y$, the $(X, Y)$-paths is a family of internally disjoint paths starting at a vertex $x\in X$, ending at a vertex $y\in Y$ and whose internal vertices belong neither to $X$ nor $Y$. If $X=\{x\}$, then the $(X, Y)$-paths is a family of internal disjoint paths whose starting vertex is $x$ and the terminal vertices are distinct in $Y$, which is referred to as a {\em $k$-fan} from $x$ to $Y$. For terminologies and notations not defined here we follow ~\cite{B}.

The hypercube is one of the most fundamental interconnection network. An $n$-dimensi\\onal hypercube is an undirected graph $Q_{n}=(V,E)$ with $|V|=2^{n}$ and $|E|=n2^{n-1}.$ Each vertex can be represented by an $n$-bit binary string. There is an edge between two vertices whenever their binary string representation differs in only one bit position. The dual-cube was introduced by Li and Peng in \cite{L}. As an invariant of the hypercube, it not only keeps numerous desirable properties of the hypercube, but also reduces the interconnection complexity.  The $n$-dimensional dual cube, denoted by $D_{n}$, has $2^{2n-1}$ vertices, each labeled by a $(2n-1)$-bit binary string $u_{1}u_{2}\cdots u_{2n-1}$, and $u_{i}\in \{0,1\}$ for $i=1,2\cdots, 2n-1.$ Two vertices $u=u_{1}u_{2}\cdots u_{2n-1}$ and $v=v_{1}v_{2}\cdots v_{2n-1}$ are adjacent if and only if the following conditions are satisfied:

$(1)$ $u$ and $v$ differ in exactly one bit position $i$.

$(2)$ If $1\leq i\leq n-1$, then $u_{2n-1}=v_{2n-1}=0$.

$(3)$ If $n\leq i\leq 2n-2$, then $u_{2n-1}=v_{2n-1}=1$.

By the definition of the $n$-dimensional dual cube, we can see that the set of vertices with fixed bits in positions $n,\cdots, 2n-2$ and last bit $0$ form an $(n-1)$-dimensional hypercube, which is called a cluster of class $0$. And the set of vertices with fixed bits in positions $1,\cdots, n-1$ and last bit $1$ also form an $(n-1)$-dimensional hypercube, which is called a cluster of class $1$. Edges connecting two vertices in different clusters of different class are called cross edges. Let $v$ be a vertex of cluster $0 (1)$, then we call the neighbour of the vertex $v$ in cluster $1 (0)$ the outside neighbour of $v$ in $D_{n}$. The $n$-dimensional dual cube has $2^{n-1}$ clusters of class $0$ and $2^{n-1}$ clusters of class $1$. For convenience, we denote the clusters of class $0$ by $D_{01}, D_{02},\cdots $ and $D_{02^{n-1}}$, respectively. And we denote the clusters of class $1$ by $D_{11}, D_{12},\cdots$ and $D_{12^{n-1}}$, respectively. The $3$-dimensional dual cube is shown as figure $1$.

\begin{figure}[!ht]
\begin{center}
%\vskip4cm
\includegraphics[scale=0.6]{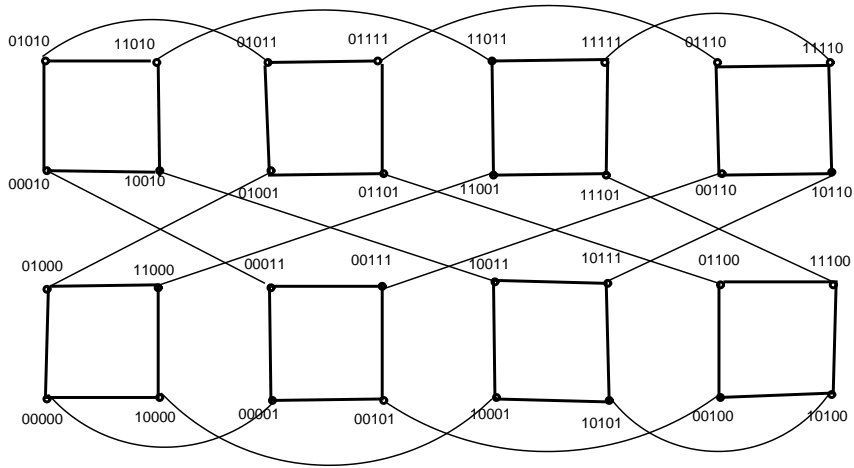}
\end{center}
\vskip0.5cm
\caption{The $3-$dimensional dual cube}\label{F1}
\end{figure}

The $n$-dimensional dual cube $D_{n}$ is an $n$-regular graph. It consists of $2^{n}$ clusters, half of them are class $0$, and the other half of them is class $1$. There are many research results about dual cubes, one can refer \cite{Ar, Yang2, Ya} etc. for the detail.

\section{The generalized $4$-connectivity  of dual cubes}

In this section, we will study the generalized $4$-connectivity  of dual cubes. The following results are useful to the main result.

\begin{lem}\label{lem1} Let $D_{n}$ be an $n$-dimensional dual cube for $n\geq 2$. Then the following results hold.
\begin{enumerate}
\item [{\rm (1)}] Any two different vertices in the same cluster of $D_{n}$ have different outside neighbours and they belong to different clusters of $D_{n}$.

\item [{\rm (2)}] There is exactly one cross edge between two clusters of different class in $D_{n}$.
\end{enumerate}
\end{lem}
\f {\bf Proof.} (1) Let $D_{01}, D_{02},\cdots, D_{02^{n-1}}$ be the clusters of class $0$ and $D_{11}, D_{12},\cdots, D_{12^{n-1}}$ be the clusters of class $1$ of $D_{n}$. Without loss of generality, let $u,v\in V(D_{01})$ and $u\neq v$. Let $u=u_{1}u_{2}\cdots u_{n-1}u_{n}\cdots u_{2n-2}0$ and $v=v_{1}v_{2}\cdots v_{n-1}v_{n}\cdots v_{2n-2}0$. As $u\neq v$, then $u_{1}u_{2}\cdots u_{n-1}$ and $v_{1}v_{2}\cdots v_{n-1}$ have at least one different digit. The outside neighbour $u^{\prime}$ of $u$ is $u^{\prime}=u_{1}u_{2}\cdots u_{n-1}u_{n}\cdots u_{2n-2}1$ and
the outside neighbour $v^{\prime}$ of $v$ is $v^{\prime}=v_{1}v_{2}\cdots v_{n-1}v_{n}\cdots v_{2n-2}1$. If $u^{\prime}$ and $v^{\prime}$ belong to the same cluster, then $u_{1}u_{2}\cdots u_{n-1}=v_{1}v_{2}\cdots v_{n-1}$ by the definition of $D_{n}$, which is a contradiction.

(2) By $(1)$, the result holds easily.
\hfill\qed
\begin{lem}{\rm(\cite{B})}\label{lem2}
 $\kappa(Q_{n})=n$ for $n\geq 2$.
\end{lem}

\begin{theorem}{\rm(\cite{L})}\label{lem3}
 The generalized $4$-connectivity of hypercube $Q_{n}$ is $n-1$, where $n\geq 2$.
\end{theorem}

\begin{lem}{\rm(\cite{L})}\label{lem4}
Let $G$ be an $r$-regular graph. If $\kappa_{k}(G)=r-1$, then $\kappa_{k-1}(G)=r-1$, where $k\geq 4$.
\end{lem}

\begin{lem}{\rm(\cite{l5})}\label{lem5}
If there are two adjacent vertices of degree $\delta(G)$, then $\kappa_{k}(G)\leq \delta-1$ for $3\leq k \leq |V(G)|$.
\end{lem}

\begin{lem}{\rm(\cite{B})}\label{lem6}
Let $G=(V, E)$ be a $k$-connected graph, let $x$ be a vertex of $G$, and let $Y\subseteq V\setminus \{x\}$ be a set of at least $k$ vertices of $G$. Then there exists a $k$-fan in $G$ from $x$ to $Y$. That is, there exists a family of $k$ internally vertex-disjoint $(x, Y)$-paths whose terminal vertices are distinct in $Y$.
\end{lem}

%\begin{lem}{\rm(\cite{Ar})}\label{lem7} $\kappa(D_{n})=n$ for $\geq 2$.
%\end{lem}
\begin{theorem}\label{lem8}
Let $H=D_{n}[(\bigcup_{j=1}^{k}V(D_{0i_{j}}))\bigcup(\bigcup_{i=1}^{l}V(D_{1j_{i}}))]$, where $k\geq 1, l\geq 1,$ $D_{0i_{j}}$ is a cluster of class $0$ of $D_{n}$ and $D_{1j_{i}}$ is a cluster of class $1$ of $D_{n}$. Then $H$ is connected for $n\geq 3$.
\end{theorem}
\f {\bf Proof.} Without loss of generality, let $H=D_{n}[(\bigcup_{j=1}^{k}V(D_{0j}))\bigcup(\bigcup_{i=1}^{l}V(D_{1i}))]$, where $k\geq 1, l\geq 1,$ $D_{0j}$ is a cluster of class $0$ of $D_{n}$ and $D_{1i}$ is a cluster of class $1$ of $D_{n}$. To prove the result, we just need to show for any two distinct
vertices $x$ and $y$ of $H$, there is a path between $x$ and $y$ in $H$.

Case $1$. $x$ and $y$ belong to the same cluster.

Without loss of generality, let $x, y\in V(D_{01})$. As $D_{01}$ is connected. Then there is a path $P$ between $x$ and $y$ in $D_{01}$. As $D_{01}$ is a subgraph of $H$, then $P$ is a path between $x$ and $y$ in $H$.

Case $2$. $x$ and $y$ belong to two clusters of the same class.

Without loss of generality, let $x\in V(D_{01})$ and $y\in V(D_{02})$. By Lemma~\ref{lem1}(1), there is a vertex $u$ of $D_{01}$ such that $u^{1}$, the outside neighbour of $u$, belongs to $D_{11}$. Similarly, there is a vertex $v$ of $D_{02}$ such that $v^{1}$, the outside neighbour of $v$, belongs to $D_{11}$. As $D_{01}, D_{02}$ and $D_{11}$ are connected. Then there is a path $P$ between $x$ and $u$ in $D_{01}$, a path $Q$ between $u^{1}$ and $v^{1}$ in $D_{11}$ and a path $R$ between $y$ and $v$ in $D_{02}$. If $u=x$, then the path $P$ is the single vertex $x$ and if $v=y$, then the path $R$ is the single vertex $y$. Let $\widehat{P}=xPuu^{1}Qv^{1}vRy$, then it is a path between $x$ and $y$ in $H$.

Case $3$. $x$ and $y$ belong to two clusters of different class.

Without loss of generality, let $x\in V(D_{01})$ and $y\in V(D_{11})$. By Lemma~\ref{lem1}(1), there is a vertex $u$ of $D_{01}$ such that $u^{1}$, the outside neighbour of $u$, belongs to $D_{11}$. As $D_{01}$ and $D_{11}$ are connected. Then there is a path $P$ between $x$ and $u$ in $D_{01}$ and a path $Q$ between $u^{1}$ and $y$ in $D_{11}$. If $u=x$, then the path $P$ is the single vertex $x$ and if $u^{1}=y$, then the path $Q$ is the single vertex $y$. Let $\widehat{P}=xPuu^{1}Qy$, then it is a path between $x$ and $y$ in $H$.
\hfill\qed

To prove the generalized $4$-connectivity of dual cubes, the following results are useful.

\begin{theorem}\label{Lem81}
Let $D_{n}$ be an $n$-dimensional dual cube and $S=\{x, y, z, w\}$, where $x, y, z$ and $w$ are any four distinct vertices of $D_{n}$ for $n\geq 4$. If the vertices in $S$ belong to two clusters of $D_{n}$, then there are $n-1$ internally disjoint trees connecting $S$ in $D_{n}$.
\end{theorem}
\f {\bf Proof.} Let $D_{01}, D_{02},\cdots, D_{02^{n-1}}$ be the clusters of class $0$ and $D_{11}, D_{12},\cdots$ and $D_{12^{n-1}}$ be the clusters of class $1$ of $D_{n}$. Let $S=\{x, y, z, w\}$, where $x, y, z$ and $w$ are any four distinct vertices of $D_{n}$ for $n\geq 4$. Let the vertices in $S$ belong to two clusters of $D_{n}$. By the symmetry of $D_{n}$, the following cases are considered.

Case $1$. $|S\bigcap V(D_{0i})|=3$ and $|S\bigcap V(D_{0j})|=1$ for distinct $i,j\in\{1,2,\cdots, 2^{n-1}\}$.

Without loss of generality, let $|S\bigcap V(D_{01})|=3$ and $|S\bigcap V(D_{02})|=1$. Let $\{x,y,z\}\subseteq V(D_{01})$ and $w\in V(D_{02})$. See Figure $2$. As $x\neq y$, without loss of generality, we may assume that the $(n-1)$-th digit of $x$ and $y$ are different. Hence, we may assume that $x=x_{1}x_{2}\cdots x_{n-2}0H0$ and $y=y_{1}y_{2}\cdots y_{n-2}1H0$, where $H=x_{n}x_{n+1}\cdots x_{2n-2}$. Assume the $(n-1)$-th digit of $z$ is $0$ and we divide $D_{01}$ along the $(n-1)$th dimension into two copies of $Q_{n-2}$, denoted by $Q^{n-1}[0]$ and $Q^{n-1}[1]$, respectively. Then $x, z\in V(Q^{n-1}[0])$ and $y\in V(Q^{n-1}[1])$. By Lemma~\ref{lem2}, $\kappa(Q_{n-2})=n-2$. Then there are $n-2$ internally disjoint paths $P_{1}, P_{2},\cdots, P_{n-2}$ between $x$ and $z$ in $Q^{n-1}[0]$. Let $x_{i}\in V(P_{i})$ such that $y_{i}\in V(Q^{n-1}[1])$ and $y_{i}\neq y$, where $y_{i}$ is the neighbour of $x_{i}$ in $Q^{n-1}[1]$. This can be done as any two different vertices in $Q^{n-1}[0]$ have different neighbours in $Q^{n-1}[1]$ and $P_{i}$s are internally disjoint for $1\leq i\leq n-2$. Let $Y=\{y_{1}, y_{2}, \ldots, y_{n-2}\}$. By Lemma~\ref{lem2}, $\kappa(Q^{n-1}[1])=n-2$. By Lemma~\ref{lem6}, there are $n-2$ internally disjoint paths $P_{1}^{\prime}, P_{2}^{\prime},\cdots, P_{n-2}^{\prime}$ from $y$ to $Y$. Let $\widehat{T}_{i}=P_{i}\bigcup x_{i}y_{i}\bigcup P_{i}^{\prime}$ for each $i\in\{1,2,\cdots, n-2\}$, then $n-2$ internally disjoint trees $\widehat{T}_{i}$s for $1\leq i\leq n-2$ that connecting $x, y$ and $z$ are obtained in $D_{01}$.

\begin{figure}[!ht]
\begin{center}
%\vskip4cm
\includegraphics[scale=0.7]{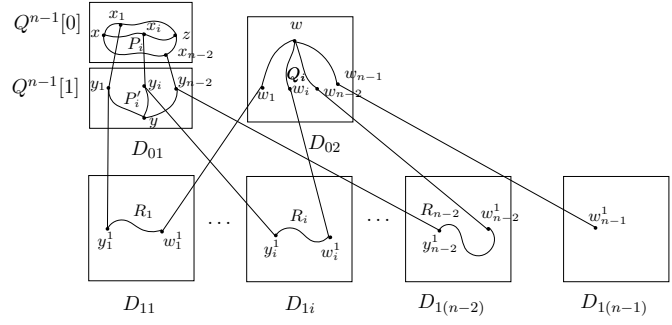}
\end{center}
\vskip0.1cm
\caption{The illustration of Case $1$}\label{F2} %for $y_{i}^{0}\neq w$ for each $i\in\{1,2,\cdots, n-2\}$
\end{figure}

Let $y_{i}^{1}$ be the outside neighbour of $y_{i}$ in $D_{n}$ for each $i\in\{1, 2, \cdots, n-2\}$. Without loss of generality, let $y_{i}^{1}\in V(D_{1i})$ for $1\leq i\leq n-2$. Then choose $n-1$ vertices $w_{1}, w_{2},\cdots, w_{n-1}$ from $D_{02}$ such that $w_{i}^{1}$, the outside neighbour of $w_{i}$ in $D_{n}$, belongs to $D_{1i}$ for each $i\in\{1, 2, \cdots, n-1\}$. For convenience, let $W=\{w_{1}, w_{2},\cdots, w_{n-1}\}$. By Lemma~\ref{lem2}, $\kappa(D_{02})=n-1$ as $D_{02}$ is a copy of $Q_{n-1}$. By Lemma~\ref{lem6}, there are $n-1$ internally disjoint paths $Q_{1}, Q_{2}, \cdots, Q_{n-1}$ from $w$ to $W$ and let the terminal vertex of $Q_{i}$ is $w_{i}$ for each $i\in\{1, 2, \cdots, n-2\}$. As $D_{1i}$ is connected for $1\leq i\leq n-2$, then there is a path $R_{i}$ between $y_{i}^{1}$ and $w_{i}^{1}$ in $D_{1i}$ for each $i\in\{1, 2, \cdots, n-2\}$.

Let $T_{i}=\widehat{T}_{i}\bigcup Q_{i}\bigcup R_{i}\bigcup y_{i}y_{i}^{1}\bigcup w_{i}w_{i}^{1}$ for each $i\in\{1, 2, \cdots, n-2\}$. Then $n-2$ internally disjoint trees $T_{i}$s for $1\leq i\leq n-2$ that connecting $x, y, z$ and $w$ are obtained. Let $x^{1}, y^{1}$ and $z^{1}$ be the outside neighbours of $x, y$ and $z$ in $D_{n}$, respectively. By Lemma~\ref{lem1}(2), $\{x^{1}, y^{1}, z^{1}, w_{n-1}^{1}\}\subset V(D_{n}-\bigcup_{i=1}^{2}D_{0i}-\bigcup_{i=1}^{n-2}D_{1i})$. As $2^{n-1}> n-2$ for $n\geq 4$. By Theorem~\ref{lem8}, $D_{n}-\bigcup_{i=1}^{2}D_{0i}-\bigcup_{i=1}^{n-2}D_{1i}$ is connected. Then there is a tree $\widetilde{T}_{n-1}$ connecting $x^{1}, y^{1}, z^{1}$ and $w_{n-1}^{1}$ in $D_{n}-\bigcup_{i=1}^{2}D_{0i}-\bigcup_{i=1}^{n-2}D_{1i}$. Let $T_{n-1}=\widetilde{T}_{n-1}\bigcup Q_{n-1}\bigcup xx^{1}\bigcup yy^{1}\bigcup zz^{1}\bigcup w_{n-1}w_{n-1}^{1}$, then $n-1$ internally disjoint $S$-trees $T_{i}$s for $1\leq i\leq n-1$ are obtained in $D_{n}$.

Case $2$. $|S\bigcap V(D_{0i})|=3$ and $|S\bigcap V(D_{1j})|=1$ for $i,j\in\{1,2,\cdots, 2^{n-1}\}$.

Without loss of generality, let $|S\bigcap V(D_{01})|=3$ and $|S\bigcap V(D_{11})|=1$. Let $\{x,y,z\}\subseteq V(D_{01})$ and $w\in V(D_{11})$. See Figure $3$. Let $Y$ and $\widehat{T}_{i}$s for $1\leq i\leq n-2$ be the same as case $1$. Note that $Y=\{y_{1}, y_{2}, \ldots, y_{n-2}\}$ and $\widehat{T}_{1}, \widehat{T}_{2}, \cdots, \widehat{T}_{n-2}$ are $n-2$ internally disjoint trees connecting $x, y$ and $z$ in $D_{01}$. Let $Y^{1}=\{y_{1}^{1}, y_{2}^{1},\ldots, y_{n-2}^{1}\}$, where $y_{i}^{1}$ is the outside neighbour of $y_{i}$ in $D_{n}$ for each $y_{i}\in Y$.

Subcase $2.1$. If no element in $Y^{1}$ belongs to $D_{11}$. Without loss of generality, let $y_{i}^{1}\in V(D_{1(i+1)})$ for each $i\in\{1,2,\cdots, n-2\}$. Let $N_{D_{11}}(w)=\{w_{1}, w_{2},\cdots,w_{n-2}, w_{n-1}\}$. Then choose $n-2$ vertices from $N_{D_{11}}(w)$ such that none of them with outside neighbour in $D_{01}$, say $w_{1}, w_{2},\cdots,w_{n-2}$ and let $W=\{w_{1}, w_{2},\cdots,w_{n-2}\}$. Without loss of generality, let $w_{i}^{0}\in V(D_{0(i+1)})$, where $w_{i}^{0}$ is the outside neighbour  of $w_{i}$ in $D_{n}$. By Theorem~\ref{lem8}, $D_{n}[V(D_{0(i+1)})\bigcup V(D_{1(i+1)})]$ is connected. Then there is a path $Q_{i}$ between $w_{i}^{0}$ and $y_{i}^{1}$ in $D_{n}[V(D_{0(i+1)})\bigcup V(D_{1(i+1)})]$ for each $i\in\{1,2,\cdots, n-2\}$. Let $T_{i}=\widehat{T}_{i}\bigcup Q_{i}\bigcup y_{i}y_{i}^{1}\bigcup w_{i}w_{i}^{0}\bigcup \\w_{i}w$ for each $i\in\{1,2,\cdots, n-2\}$. Then $n-2$ internally disjoint $S$-trees $T_{i}$s for $1\leq i\leq n-2$ that connecting $x, y, z$ and $w$ are obtained. Let $x^{1}, y^{1}, z^{1}, w^{0}$ and $w_{n-1}^{0}$ be the outside neighbours of $x, y, z, w$ and $w_{n-1}$ in $D_{n}$, respectively.

Subcase $2.1.1$. If one of $x^{1}, y^{1}$ and $z^{1}$ belongs to $D_{11}$, say $x^{1}$. By Lemma~\ref{lem2}, $\kappa(D_{11})=n-1$. Then $D_{11}\setminus W$ is connected. Then there is a path $Q$ between $x^{1}$ and $w$ in $D_{11}\setminus W$.

If $x^{1}\neq w$. By Lemma~\ref{lem1}(2), $\{y^{1}, z^{1}, w^{0}\}\subseteq V(D_{n}-\bigcup_{i=1}^{n-1}D_{0i}-\bigcup_{i=1}^{n-1}D_{1i})$. As $2^{n-1}> n-1$ for $n\geq 4$. By Theorem~\ref{lem8}, $D_{n}-\bigcup_{i=1}^{n-1}D_{0i}-\bigcup_{i=1}^{n-1}D_{1i}$ is connected. Then there is a tree $\widetilde{T}_{n-1}$ connecting $y^{1}, z^{1}$ and $w^{0}$ in $D_{n}-\bigcup_{i=1}^{n-1}D_{0i}-\bigcup_{i=1}^{n-1}D_{1i}$. Let $T_{n-1}= xx^{1}\bigcup Q\bigcup ww^{0}\bigcup yy^{1}\bigcup zz^{1}\bigcup \widetilde{T}_{n-1}$, then $n-1$ internally disjoint $S$-trees $T_{i}$s for $1\leq i\leq n-1$ are obtained in $D_{n}$.

If $x^{1}=w$. By Lemma~\ref{lem1}(2), $\{w_{n-1}^{0}, y^{1}, z^{1}\}\subseteq V(D_{n}-\bigcup_{i=1}^{n-1}D_{0i}-\bigcup_{i=1}^{n-1}D_{1i})$. As $2^{n-1}> n-1$ for $n\geq 4$. By Theorem~\ref{lem8}, $D_{n}-\bigcup_{i=1}^{n-1}D_{0i}-\bigcup_{i=1}^{n-1}D_{1i}$ is connected. Then there is a tree $\widetilde{T}_{n-1}$ connecting $y^{1}, z^{1}$ and $w_{n-1}^{0}$ in $D_{n}-\bigcup_{i=1}^{n-1}D_{0i}-\bigcup_{i=1}^{n-1}D_{1i}$. Let $T_{n-1}= xw\bigcup ww_{n-1}\bigcup w_{n-1}w_{n-1}^{0}\bigcup \widetilde{T}_{n-1}\bigcup yy^{1}\bigcup zz^{1}$. Then $n-1$ internally disjoint $S$-trees $T_{i}$s for $1\leq i\leq n-1$ are obtained in $D_{n}$.

Subcase $2.1.2$. If none of $x^{1}, y^{1}$ and $z^{1}$ belongs to $D_{11}$. Then $\{x^{1}, y^{1}, z^{1}\}\subseteq V(D_{n}-\bigcup_{i=1}^{n-1}D_{0i}-\bigcup_{i=1}^{n-1}D_{1i})$. By Lemma~\ref{lem1}(2), at least one of $w^{0}$ and $w_{n-1}^{0}$ does not belong to $D_{01}$.

If $w^{0}\notin V(D_{01})$, then $w^{0}\in V(D_{n}-\bigcup_{i=1}^{n-1}D_{0i}-\bigcup_{i=1}^{n-1}D_{1i})$. By Theorem~\ref{lem8}, $D_{n}-\bigcup_{i=1}^{n-1}D_{0i}-\bigcup_{i=1}^{n-1}D_{1i}$ is connected. Then there is a tree $\widetilde{T}_{n-1}$ connecting $x^{1}, y^{1}$ and $z^{1}$ and $w^{0}$ in $D_{n}-\bigcup_{i=1}^{n-1}D_{0i}-\bigcup_{i=1}^{n-1}D_{1i}$. Let $T_{n-1}= \widetilde{T}_{n-1}\bigcup xx^{1}\bigcup yy^{1}\bigcup zz^{1}\bigcup ww^{0}$, then $n-1$ internally disjoint $S$-trees $T_{i}$s for $1\leq i\leq n-1$ are obtained in $D_{n}$.

If $w_{n-1}^{0}\notin V(D_{01})$, then $w_{n-1}^{0}\in V(D_{n}-\bigcup_{i=1}^{n-1}D_{0i}-\bigcup_{i=1}^{n-1}D_{1i})$. By Theorem~\ref{lem8}, $D_{n}-\bigcup_{i=1}^{n-1}D_{0i}-\bigcup_{i=1}^{n-1}D_{1i}$ is connected. Then there is a tree $\widetilde{T}_{n-1}$ connecting $x^{1}, y^{1}$ and $z^{1}$ and $w_{n-1}^{0}$ in $D_{n}-\bigcup_{i=1}^{n-1}D_{0i}-\bigcup_{i=1}^{n-1}D_{1i}$. Let $T_{n-1}= \widetilde{T}_{n-1}\bigcup xx^{1}\bigcup yy^{1}\bigcup zz^{1}\bigcup w_{n-1}w_{n-1}^{0} \bigcup ww_{n-1}$, then $n-1$ internally disjoint $S$-trees $T_{i}$s for $1\leq i\leq n-1$ are obtained in $D_{n}$.

Subcase $2.2$. If there is some element in $Y^{1}$ belongs to $D_{11}$. Without loss of generality, let $y_{1}^{1}\in V(D_{11})$ and $y_{i}^{1}\in V(D_{1i})$ for $2\leq i\leq n-2$. Define $W$ same as Subcase $2.1$. By Lemma~\ref{lem2}, $\kappa(D_{11})=n-1$. Then $D_{11}\setminus W$ is connected. Then there is a path $Q$ between $y_{1}^{1}$ and $w$ in $D_{11}\setminus W$. Then let $T_{1}=\widehat{T}_{1}\bigcup y_{1}y_{1}^{1}\bigcup Q$. By Theorem~\ref{lem8}, $D_{n}[V(D_{0(i+1)})\bigcup V(D_{1i})]$ is connected. Then there is a path $Q_{i}$ between $w_{i}^{0}$ and $y_{i}^{1}$ in $D_{n}[V(D_{0(i+1)})\bigcup V(D_{1i})]$ for each $i\in\{2,\cdots, n-2\}$. Let $T_{i}=\widehat{T}_{i}\bigcup Q_{i}\bigcup y_{i}y_{i}^{1}\bigcup w_{i}w_{i}^{0}\bigcup w_{i}w$ for each $2\leq i\leq n-2$. Combining with $T_{1}$, then $n-2$ internally disjoint $S$-trees $T_{i}$s for $1\leq i\leq n-2$ that connecting $x, y, z$ and $w$ are obtained.

Subcase $2.2.1$. If $y_{1}^{1}\neq w$. By Lemma~\ref{lem1}(2), $\{x^{1}, y^{1}, z^{1}, w^{1} \}\subseteq V(D_{n}-\bigcup_{i=1}^{n-1}D_{0i}-\bigcup_{i=1}^{n-1}D_{1i})$. As $2^{n-1}> n-1$ for $n\geq 4$. By Theorem~\ref{lem8}, $D_{n}-\bigcup_{i=1}^{n-1}D_{0i}-\bigcup_{i=1}^{n-1}D_{1i}$ is connected. Then there is a tree $\widetilde{T}_{n-1}$ connecting $x^{1}, y^{1}$ and $z^{1}$ and $w^{1}$ in $D_{n}-\bigcup_{i=1}^{n-1}D_{0i}-\bigcup_{i=1}^{n-1}D_{1i}$. Let $T_{n-1}= \widetilde{T}_{n-1}\bigcup xx^{1}\bigcup yy^{1}\bigcup zz^{1}\bigcup ww^{1}$, then $n-1$ internally disjoint $S$-trees $T_{i}$s for $1\leq i\leq n-1$ are obtained in $D_{n}$.

Subcase $2.2.2$. If $y_{1}^{1}= w$. Then $T_{1}=\widehat{T}_{1}\bigcup y_{1}y_{1}^{1}$. By Lemma~\ref{lem1}(2), $\{w_{n-1}^{0}, x^{1}, y^{1}, z^{1}\}\subseteq V(D_{n}-\bigcup_{i=1}^{n-1}D_{0i}-\bigcup_{i=1}^{n-1}D_{1i})$. As $2^{n-1}> n-1$ for $n\geq 4$. By Theorem~\ref{lem8}, $D_{n}-\bigcup_{i=1}^{n-1}D_{0i}-\bigcup_{i=1}^{n-1}D_{1i}$ is connected. Then there is a tree $\widetilde{T}_{n-1}$ connecting $w_{n-1}^{0}, x^{1}, y^{1}$ and $z^{1}$ in $D_{n}-\bigcup_{i=1}^{n-1}D_{0i}-\bigcup_{i=1}^{n-1}D_{1i}$. Let $T_{n-1}= \widetilde{T}_{n-1}\bigcup xx^{1}\bigcup yy^{1}\bigcup zz^{1}\bigcup ww_{n-1}\bigcup w_{n-1}w_{n-1}^{0}$, then $n-1$ internally disjoint $S$-trees $T_{i}$s for $1\leq i\leq n-1$ are obtained in $D_{n}$.

\begin{figure}[!ht]
\begin{center}
%\vskip4cm
\includegraphics[scale=0.7]{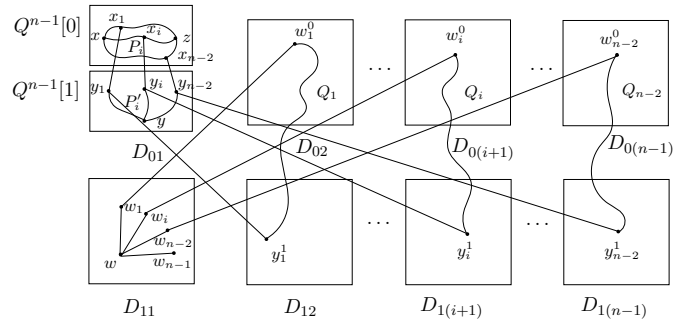}
\end{center}
\vskip0.1cm
\caption{The illustration of Subcase $2.1$}\label{F3}% for $w_{i}\neq w$ for each $i\in\{1,2,\cdots, n-2\}$
\end{figure}

Case $3$. $|S\bigcap V(D_{0i})|=2$ and $|S\bigcap V(D_{0j})|=2$ for $i\neq j$ and $i,j\in\{1,2,\cdots, 2^{n-1}\}$.

Without loss of generality, let $|S\bigcap V(D_{01})|=2$ and $|S\bigcap V(D_{02})|=2$. Let $\{x,y\}\subseteq V(D_{01})$ and $\{z, w\}\subseteq V(D_{02})$. Recall that $D_{ij}$ is a copy of $Q_{n-1}$. By Lemma~\ref{lem2}, $\kappa(D_{01})=\kappa(D_{02})=n-1$. Then there are $n-1$ internally disjoint paths $P_{1}, P_{2}, \cdots, P_{n-1} $ between $x$ and $y$ in $D_{01}$ and $n-1$ internally disjoint paths $P_{1}^{\prime}, P_{2}^{\prime}, \cdots, P_{n-1}^{\prime} $ between $z$ and $w$ in $D_{02}$. Let $x_{i}\in V(P_{i})$ and $y_{i}\in V(P_{i}^{\prime})$ for each $i\in\{1,2,\cdots, n-1\}$. Let $x_{i}^{1}$ and $y_{i}^{1}$ be the outside neighbour of $x_{i}$ and $y_{i}$ in $D_{n}$ for each $i\in\{1,2,\cdots, n-1\}$, respectively. Without loss of generality, let $x_{i}^{1}\in V(D_{1i})$ for each $i\in\{1,2,\cdots, n-1\}$. Let $X^{\prime}=\{x_{1}^{1}, x_{2}^{1},\cdots, x_{n-1}^{1}\}$ and $Y^{\prime}=\{y_{1}^{1}, y_{2}^{1},\cdots, y_{n-1}^{1}\}$, then $X^{\prime}\bigcap Y^{\prime}=\emptyset$.

Subcase $3.1$. If there is no vertex in $X^{\prime}$ belongs to the same cluster with the vertices in $Y^{\prime}$. Then without loss of generality, let $y_{i}^{1}\in V(D_{1(i+n-1)})$ for each $i\in\{1,2,\cdots, n-1\}$. This can be done as $2^{n-1}> 2(n-1)$ for $n\geq 4$. See Figure $4$. By Theorem~\ref{lem8}, $D_{n}[V(D_{1i})\bigcup V(D_{1(i+n-1)})\bigcup V(D_{0(i+n-1)})]$ is connected for each $i\in\{1,2,\cdots, n-1\}$. Then there is a path $Q_{i}$ between $x_{i}^{1}$ and $y_{i}^{1}$ in $D_{n}[V(D_{1i})\bigcup V(D_{1(i+n-1)})\bigcup V(D_{0(i+n-1)})]$ for each $i\in\{1,2,\cdots, n-1\}$. Let $T_{i}=P_{i}\bigcup P_{i}^{\prime}\bigcup Q_{i}\bigcup x_{i}x_{i}^{1}\bigcup y_{i}y_{i}^{1}$ for each $i\in\{1,2,\cdots, n-1\}$. Then $n-1$ internally disjoint $S$-trees $T_{i}$s for $1\leq i\leq n-1$ are obtained in $D_{n}$.

Subcase $3.2$. If there are some vertices in $X^{\prime}$ belong to the same cluster with some vertices in $Y^{\prime}$.
Then without loss of generality, let $x_{i}^{1}$ and $y_{i}^{1}$ belong to the same cluster for $t\leq i\leq n-1$ and let $y_{i}^{1}\in V(D_{1(i+n-1)})$ for $1\leq i\leq t-1$, where $1\leq t\leq n-1$. By Theorem~\ref{lem8}, $D_{n}[V(D_{1i})\bigcup V(D_{1(i+n-1)})\bigcup V(D_{0(i+n-1)})]$ is connected for $1\leq i\leq t-1$. Then there is a path $Q_{i}$ between $x_{i}^{1}$ and $y_{i}^{1}$ in $D_{n}[V(D_{1i})\bigcup V(D_{1(i+n-1)})\bigcup V(D_{0(i+n-1)})]$ for each $i\in\{1,2,\cdots, t-1\}$. As $x_{i}, y_{i}\in V(D_{1i})$ and $D_{1i}$ is connected for each $i\in\{t, t+1,\cdots, n-1\}$, then there is a path $R_{i}$ between $x_{i}$ and $y_{i}$ in $D_{1i}$ for $t\leq i \leq n-1$. Let $T_{i}=P_{i}\bigcup P_{i}^{\prime}\bigcup Q_{i}\bigcup x_{i}x_{i}^{1}\bigcup y_{i}y_{i}^{1}$ for $1\leq i \leq t-1$ and $T_{i}=P_{i}\bigcup P_{i}^{\prime}\bigcup R_{i}\bigcup x_{i}x_{i}^{1}\bigcup y_{i}y_{i}^{1}$ for $t\leq i \leq n-1$. Then $n-1$ internally disjoint trees $T_{i}$s for $1\leq i\leq n-1$ are obtained in $D_{n}$.
\begin{figure}[!ht]
\begin{center}
%\vskip4cm
\includegraphics[scale=0.7]{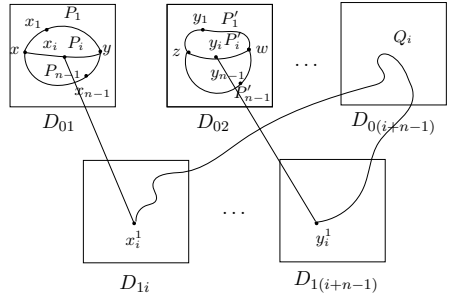}
\end{center}
\vskip0.1cm
\caption{The illustration of Subcase $3.1$}\label{F4}% when there is no vertex in $X^{\prime}$ belongs to the same cluster with the vertices in $Y^{\prime}$
\end{figure}

Case $4$. $|S\bigcap V(D_{0i})|=2$ and $|S\bigcap V(D_{1j})|=2$ for $i,j\in\{1,2,\cdots, 2^{n-1}\}$.

Without loss of generality, let $|S\bigcap V(D_{01})|=2$ and $|S\bigcap V(D_{11})|=2$. Let $\{x,y\}\subseteq V(D_{01})$ and $\{z, w\}\subseteq V(D_{11})$. See Figure $5$. Recall that $D_{ij}$ is a copy of $Q_{n-1}$. By Lemma~\ref{lem2}, $\kappa(D_{01})=\kappa(D_{11})=n-1$. Then there are $n-1$ internally disjoint paths $P_{1}, P_{2}, \cdots, P_{n-1} $ between $x$ and $y$ in $D_{01}$ and $n-1$ internally disjoint paths $P_{1}^{\prime}, P_{2}^{\prime}, \cdots, P_{n-1}^{\prime} $ between $z$ and $w$ in $D_{11}$. Let $x_{i}\in V(P_{i})$ such that $x_{i}^{1}$, the outside neighbour of $x_{i}$, does not belong to $D_{11}$. Without loss of generality, let $x_{i}^{1}\in V(D_{1(i+1)})$. Similarly, let $y_{i}\in V(P_{i}^{\prime})$ such that $y_{i}^{0}$, the outside neighbour of $y_{i}$, does not belong to $D_{01}$. Without loss of generality, let $y_{i}^{0}\in V(D_{0(i+1)})$. This can be done as $2^{n-1}> n$ for $n\geq 4$. By Theorem~\ref{lem8}, $D_{n}[V(D_{0(i+1)})\bigcup V(D_{1(i+1)})]$ is connected. Then there is a path $Q_{i}$ between $x_{i}^{1}$ and $y_{i}^{0}$ in $D_{n}[V(D_{0(i+1)})\bigcup V(D_{1(i+1)})]$ for each $i\in\{1,2,\cdots, n-1\}$. Let $T_{i}=P_{i}\bigcup P_{i}^{\prime}\bigcup Q_{i}\bigcup x_{i}x_{i}^{1}\bigcup y_{i}y_{i}^{0}$ for each $i\in\{1,2,\cdots, n-1\}$. Then $n-1$ internally disjoint $S$-trees $T_{i}$s for $1\leq i\leq n-1$ are obtained in $D_{n}$.

\begin{figure}[!ht]
\begin{center}
%\vskip4cm
\includegraphics[scale=0.7]{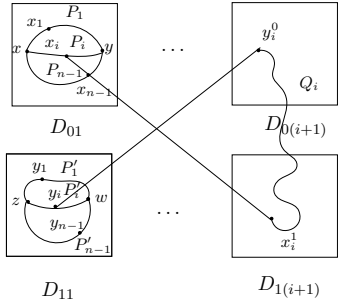}
\end{center}
\vskip0.5cm
\caption{The illustration of Case $4$}\label{F5}
\end{figure}
\hfill\qed

\begin{theorem}\label{Lem82}
Let $D_{n}$ be an $n$-dimensional dual cube and $S=\{x, y, z, w\}$, where $x, y, z$ and $w$ are any four distinct vertices of $D_{n}$ for $n\geq 4$. If the vertices in $S$ belong to three clusters of $D_{n}$, then there are $n-1$ internally disjoint trees connecting $S$ in $D_{n}$.
\end{theorem}
\f {\bf Proof.} Let $D_{01}, D_{02},\cdots, D_{02^{n-1}}$ be the clusters of class $0$ and $D_{11}, D_{12},\cdots$ and $D_{12^{n-1}}$ be the clusters of class $1$ of $D_{n}$. Let $S=\{x, y, z, w\}$, where $x, y, z$ and $w$ are any four distinct vertices of $D_{n}$ for $n\geq 4$. Let the vertices in $S$ belong to three clusters of $D_{n}$. By the symmetry of $D_{n}$, the following cases are considered.

Case $1$. $|S\bigcap V(D_{0i})|=2$, $|S\bigcap V(D_{0j})|=1$ and $|S\bigcap V(D_{0k})|=1$ for distinct $i,j,k\in\{1,2,\cdots, 2^{n-1}\}$.

Without loss of generality, let $|S\bigcap V(D_{01})|=2, |S\bigcap V(D_{02})|=1$ and $|S\bigcap V(D_{03})|=1$. Let $x, y\in V(D_{01}), z\in V(D_{02})$ and $w\in V(D_{03})$. See Figure $6$. Recall that $D_{01}$ is a copy of $Q_{n-1}$. By Lemma~\ref{lem2}, $\kappa(D_{01})=n-1$. Then there are $n-1$ internally disjoint paths $P_{1}, P_{2}, \cdots, P_{n-1} $ between $x$ and $y$ in $D_{01}$. Let $x_{i}\in V(P_{i})$ and $x_{i}^{1}$ be the outside neighbour of $x_{i}$ in $D_{n}$ for $1\leq i\leq n-1$. Without loss of generality, let $x_{i}^{1}\in V(D_{1i})$. Then choose $n-1$ vertices $z_{1}, z_{2},\cdots, z_{n-1}$ from $D_{02}$ such that $z_{i}^{1}$, the outside neighbour  of $z_{i}$, belongs to $D_{1i}$. Similarly, choose $n-1$ vertices $w_{1}, w_{2},\cdots, w_{n-1}$ from $D_{03}$ such that $w_{i}^{1}$, the outside neighbour  of $w_{i}$, belongs to $D_{1i}$ for each $i\in\{1,2,\cdots, n-1\}$. Let $Z=\{z_{1}, z_{2},\cdots, z_{n-1}\}$ and $W=\{w_{1}, w_{2},\cdots, w_{n-1}\}$. By Lemma~\ref{lem2}, $\kappa(D_{02})=\kappa(D_{03})=n-1$. By Lemma~\ref{lem6}, there are $n-1$ internally disjoint paths $Q_{1}, Q_{2},\cdots, Q_{n-1}$ from $z$ to $Z$ in $D_{02}$ and $n-1$ internally disjoint paths $R_{1}, R_{2},\cdots, R_{n-1}$ from $w$ to $W$ in $D_{03}$. It is possible that one of the paths $Q_{i}s (R_{i}s)$ is a single vertex. Let $X^{\prime}=\{x_{i}^{1}, z_{i}^{1}, w_{i}^{1}\}$, then $X^{\prime}\subset V(D_{1i})$. As $D_{1i}$ is connected, then there is a $X^{\prime}$-tree $\widehat{T}_{i}$ in $D_{1i}$ for each $i\in \{1, 2, \cdots, n-1\}$. Let $T_{i}=P_{i}\bigcup Q_{i}\bigcup R_{i}\bigcup \widehat{T}_{i} \bigcup x_{i}x_{i}^{1}\bigcup z_{i}z_{i}^{1}\bigcup w_{i}w_{i}^{1}$ for each $i\in \{1, 2, \cdots, n-1\}$. Then $n-1$ internally disjoint $S$-trees $T_{i}$s for $1\leq i\leq n-1$ are obtained in $D_{n}$.

\begin{figure}[!ht]
\begin{center}
%\vskip4cm
\includegraphics[scale=0.7]{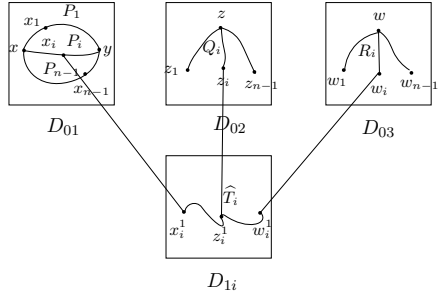}
\end{center}
\vskip0.1cm
\caption{The illustration of Case $1$}\label{F6}
\end{figure}

Case $2$. $|S\bigcap V(D_{0i})|=2$, $|S\bigcap V(D_{0j})|=1$ and $|S\bigcap V(D_{1k})|=1$ for $i\neq j$ and $i,j,k\in\{1,2,\cdots, 2^{n-1}\}$.

Without loss of generality, let $|S\bigcap V(D_{01})|=2, |S\bigcap V(D_{02})|=1$ and $|S\bigcap V(D_{11})|=1$. Let $\{x,y\}\subseteq V(D_{01}), z\in V(D_{02})$ and $w\in V(D_{11})$. Recall that $D_{01}$ is a copy of $Q_{n-1}$. By Lemma~\ref{lem2}, $\kappa(D_{01})=n-1$. Then there are $n-1$ internally disjoint paths $P_{1}, P_{2}, \cdots, P_{n-1} $ between $x$ and $y$ in $D_{01}$. See Figure $7$. Let $x_{i}\in V(P_{i})$ and $x_{i}^{1}$ be the outside neighbour of $x_{i}$ in $D_{n}$ for each $i\in \{1,2,\cdots, n-1\}$. Without loss of generality, let $x_{i}^{1}\in V(D_{1(i+1)})$. Choose $n-1$ vertices $z_{1}, z_{2},\cdots, z_{n-1}$ from $D_{02}$ such that $z_{i}^{1}$, the outside neighbour of $z_{i}$, belongs to $D_{1(i+1)}$ for each $i\in \{1,2,\cdots, n-1\}$.
Then choose $n-1$ vertices $w_{1}, w_{2},\cdots, w_{n-1}$ from $D_{11}$ such that $w_{i}^{0}$, the outside neighbour of $w$, belongs to neither $D_{01}$ nor $D_{02}$. Without loss of generality, let $w_{i}^{0}\in V(D_{0(i+2)})$. This can be done as $2^{n-1}> n+1$ for $n\geq 4$. Let $W=\{w_{1}, w_{2},\cdots, w_{n-1}\}$ and $Z=\{z_{1}, z_{2},\cdots, z_{n-1}\}$.
By Lemma~\ref{lem6}, there are $n-1$ internally disjoint paths $Q_{1}, Q_{2},\cdots, Q_{n-1}$ from $w$ to $W$ in $D_{11}$ and $n-1$ internally disjoint paths $R_{1}, R_{2},\cdots, R_{n-1}$ from $z$ to $Z$ in $D_{02}$. It is possible that one of the paths $Q_{i}s (Q_{i}s)$ is a single vertex. By Theorem~\ref{lem8}, $D_{n}[V(D_{1(i+1)})\bigcup V(D_{0(i+2)})]$ is connected. Then there is a tree $\widehat{T}_{i}$ connecting  $x_{i}^{1}, z_{i}^{1}$ and $w_{i}^{0}$ in $D_{n}[V(D_{1(i+1)})\bigcup V(D_{0(i+2)})]$. Let $T_{i}=P_{i}\bigcup Q_{i}\bigcup R_{i} \bigcup\widehat{T}_{i}\bigcup x_{i}x_{i}^{1}\bigcup z_{i}z_{i}^{1}\bigcup w_{i}w_{i}^{0}$ for each $i\in \{1, 2, \cdots, n-1\}$. Then $n-1$ internally disjoint $S$-trees $T_{i}$s for $1\leq i\leq n-1$ are obtained in $D_{n}$.

\begin{figure}[!ht]
\begin{center}
%\vskip4cm
\includegraphics[scale=0.7]{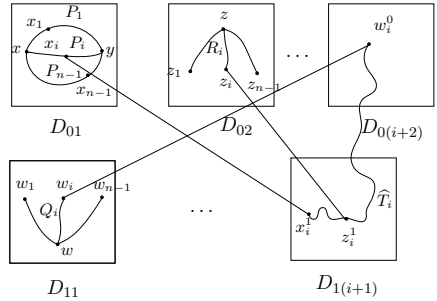}
\end{center}
\vskip0.1cm
\caption{The illustration of Case $2$}\label{F7}
\end{figure}

Case $3$. $|S\bigcap V(D_{0i})|=2$, $|S\bigcap V(D_{1j})|=1$ and $|S\bigcap V(D_{1k})|=1$ for $j\neq k$ and $i,j,k\in\{1,2,\cdots, 2^{n-1}\}$.

Without loss of generality, let $|S\bigcap V(D_{01})|=2, |S\bigcap V(D_{11})|=1$ and $|S\bigcap V(D_{12})|=1$. Let $\{x,y\}\subseteq V(D_{01}), z\in V(D_{11})$ and $w\in V(D_{12})$. Recall that $D_{01}$ is isomorphic to $Q_{n-1}$, then there are $n-1$ internally disjoint paths $P_{1}, P_{2}, \cdots, P_{n-1} $ between $x$ and $y$ in $D_{01}$. See Figure $8$.

Subcase $3.1$. If there exists $x_{i}\in V(P_{i})$ such that $x_{i}^{1}$, the outside neighbour of $x_{i}$, belongs to neither $D_{11}$ nor $D_{12}$ for each path $P_{i}$ with $i\in\{1,2,\cdots, n-1\}$. Without loss of generality, let $x_{i}^{1}\in (V(D_{1(i+2)})$. This can be done as $2^{n-1}> n+1$ for $n\geq 4$. Now, choose $n-1$ vertices $z_{1}, z_{2},\cdots, z_{n-1}$ from $D_{11}\setminus \{z\}$ such that $z_{i}^{0}$, the outside neighbour of $z_{i}$, does not belong to $D_{01}$. Without loss of generality, let $z_{i}^{0}\in V(D_{0(i+1)})$. Then choose $n-1$ vertices $w_{1}, w_{2},\cdots, w_{n-1}$ from $D_{12}\setminus \{w\}$ such that $w_{i}^{0}$, the outside neighbour of $w_{i}$, belongs to $D_{0(i+1)}$ for each $i\in\{1,2,\cdots, n-1\}$. Let $Z=\{z_{1}, z_{2},\cdots, z_{n-1}\}$ and $W=\{w_{1}, w_{2},\cdots, w_{n-1}\}$. By Lemma~\ref{lem2}, $\kappa(D_{11})=\kappa(D_{12})=n-1$. By Lemma~\ref{lem6}, there are $n-1$ internally disjoint paths $Q_{1}, Q_{2},\cdots, Q_{n-1}$ from $z$ to $Z$ in $D_{11}$ and $n-1$ internally disjoint paths $R_{1}, R_{2},\cdots, R_{n-1}$ from $w$ to $W$ in $D_{12}$. By Theorem~\ref{lem8}, $D_{n}[V(D_{0(i+1)})\bigcup V(D_{1(i+2)})]$ is connected for each $i\in\{1,2,\cdots,n-1\}$. Then there is a tree $\widehat{T}_{i}$ connecting $z_{i}^{0}, w_{i}^{0}$ and $x_{i}^{1}$ in $D_{n}[V(D_{0(i+1)})\bigcup V(D_{1(i+2)})]$. Let $T_{i}=P_{i}\bigcup Q_{i}\bigcup R_{i}\bigcup \widehat{T}_{i} \bigcup x_{i}x_{i}^{1}\bigcup z_{i}z_{i}^{0}\bigcup w_{i}w_{i}^{0}$ for each $i\in \{1, 2,\cdots, n-1\}$. Then $n-1$ internally disjoint $S$-trees $T_{i}$s for $1\leq i\leq n-1$ are obtained in $D_{n}$.

\begin{figure}[!ht]
\begin{center}
%\vskip4cm
\includegraphics[scale=0.7]{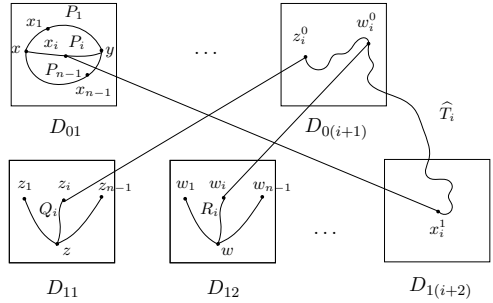}
\end{center}
\vskip0.5cm
\caption{The illustration of Subcase $3.1$}\label{F8}
\end{figure}

Subcase $3.2$. If not, then there is a path $P_{i}=xy$ with length one such that $x^{1}\in V(D_{11}), y^{1}\in V(D_{12})$ or $x^{1}\in V(D_{12}), y^{1}\in V(D_{11})$, where $x^{1}$ and $y^{1}$ are the outside neighbour of $x$ and $y$, respectively. %As $P_{i}$s are internally disjoint for $1\leq i\leq n-1$. Then there is at most one path with length one.
Without loss of generality, let $P_{n-1}=xy, x^{1}\in V(D_{11})$ and $y^{1}\in V(D_{12})$. Let $x_{i}\in V(P_{i})$ and $x_{i}^{1}$ be the outside neighbour of $x_{i}$ for $1\leq i\leq n-2$. Without loss of generality, let $x_{i}^{1}\in (V(D_{1(i+2)})$ for $1\leq i\leq n-2$.
Let $N_{D_{11}}(z)=\{z_{1}, z_{2},\cdots, z_{n-2}, z_{n-1}\}$, $N_{D_{12}}(w)=\{w_{1}, w_{2},\cdots, w_{n-2}, w_{n-1}\}$, $Z=\{z_{1}, z_{2},\cdots, z_{n-2}\}$ and $W=\{w_{1}, w_{2},\cdots, w_{n-2}\}$. By Lemma~\ref{lem2}, $\kappa(D_{11})=\kappa(D_{12})=n-1$. Then $D_{11}\setminus Z$ and $D_{12}\setminus W$ are both connected. Then there is a path $Q$ between $x^{1}$ and $z$ in $D_{11}\setminus Z$ and a path $R$ between $y^{1}$ and $w$ in $D_{12}\setminus W$. Let $T_{1}=P_{n-1}\bigcup xx^{1}\bigcup yy^{1}\bigcup Q\bigcup R$. Let $z_{i}^{0}$ and $w_{i}^{0}$ be the outside neighbour of $z_{i}$ and $w_{i}$, respectively. Without loss of generality, let $z_{i}^{0}\in V(D_{0(i+1)})$. Let $Z^{0}=\{z_{1}^{0}, z_{2}^{0},\cdots, z_{n-2}^{0}\}$ and $W^{0}=\{w_{1}^{0}, w_{2}^{0},\cdots, w_{n-2}^{0}\}$. Then $Z^{0}\bigcap W^{0}=\emptyset$.

If there is no vertex in $Z^{0}$ belongs to the same cluster with the vertices in $W^{0}$. Then without loss of generality, let $w_{i}^{0}\in V(D_{0(i+n-1)})$ for each $i\in\{1,2,\cdots, n-2\}$. This can be done as $2^{n-1}> 2(n-1)$ for $n\geq 4$. By Theorem~\ref{lem8}, $D_{n}[V(D_{1(i+2)})\bigcup V(D_{0(i+1)}\bigcup V(D_{0(i+n-1)}]$ is connected for each $i\in\{1,2,\cdots, n-2\}$. Then there is a tree $\widehat{T}_{i}$ connecting $z_{i}^{0}, w_{i}^{0}$ and $x_{i}^{1}$ in $D_{n}[V(D_{1(i+2)})\bigcup V(D_{0(i+1)}\bigcup V(D_{0(i+n-1)}]$ for each $i\in\{1,2,\cdots, n-2\}$. Let $T_{i}=P_{i}\bigcup \widehat{T}_{i}\bigcup x_{i}x_{i}^{1}\bigcup z_{i}z_{i}^{0}\bigcup w_{i}w_{i}^{0}\bigcup zz_{i}\bigcup ww_{i}$ for each $i\in\{1,2,\cdots, n-2\}$. Combining with $T_{1}$, then $n-1$ internally disjoint $S$-trees $T_{i}$s for $1\leq i\leq n-1$ are obtained in $D_{n}$.

If there are some vertices in $Z^{0}$ belong to the same cluster with some vertices in $W^{0}$.
Then without loss of generality, let $z_{i}^{0}$ and $w_{i}^{0}$ belong to the same cluster for $t\leq i\leq n-2$ and let $w_{i}^{0}\in V(D_{0(i+n-1)})$ for $1\leq i\leq t-1$, where $1\leq t\leq n-2$. By Theorem~\ref{lem8}, $D_{n}[V(D_{1(i+2)})\bigcup V(D_{0(i+1)})\bigcup V(D_{0(i+n-1)})]$ is connected for $1\leq i\leq t-1$. Then there is a tree $\widehat{T}_{i}$ connecting $z_{i}^{0}, w_{i}^{0}$ and $x_{i}^{1}$ in $D_{n}[V(D_{1(i+2)})\bigcup V(D_{0(i+1)})\bigcup V(D_{0(i+n-1)})]$ for each $i\in\{1,2,\cdots, t-1\}$. As $z_{i}^{0}, w_{i}^{0}\in V(D_{0(i+1)})$ and $D_{n}[V(D_{1(i+2)})\bigcup V(D_{0(i+1)})]$ is connected for each $i\in\{t, t+1,\cdots, n-1\}$, then there is a tree $\widehat{T}_{i}$ connecting $z_{i}^{0}, w_{i}^{0}$ and $x_{i}^{1}$ in $D_{n}[V(D_{1(i+2)})\bigcup V(D_{0(i+1)})]$ for each $t\leq i\leq n-2$. Let $T_{i}=P_{i}\bigcup \widehat{T}_{i}\bigcup x_{i}x_{i}^{1}\bigcup z_{i}z_{i}^{0}\bigcup zz_{i}\bigcup \\ww_{i}\bigcup w_{i}w_{i}^{0}$ for each $i\in\{1,2,\cdots, n-2\}$. Combining with $T_{1}$, then $n-1$ internally disjoint $S$-trees $T_{i}$s for $1\leq i\leq n-1$ are obtained in $D_{n}$.
\hfill\qed

\begin{theorem}\label{Lem83}
Let $D_{n}$ be an $n$-dimensional dual cube and $S=\{x, y, z, w\}$, where $x, y, z$ and $w$ are any four distinct vertices of $D_{n}$ for $n\geq 4$. If the vertices in $S$ belong to four clusters of $D_{n}$, then there are $n-1$ internally disjoint trees connecting $S$ in $D_{n}$.
\end{theorem}
\f {\bf Proof.} Let $D_{01}, D_{02},\cdots, D_{02^{n-1}}$ be the clusters of class $0$ and $D_{11}, D_{12},\cdots$ and $D_{12^{n-1}}$ be the clusters of class $1$ of $D_{n}$. Let $S=\{x, y, z, w\}$, where $x, y, z$ and $w$ are any four distinct vertices of $D_{n}$ for $n\geq 4$. Let the vertices in $S$ belong to four clusters of $D_{n}$. By the symmetry of $D_{n}$, the following cases are considered.

Case $1$. $|S\bigcap V(D_{0i})|=1, |S\bigcap V(D_{0j})|=1, |S\bigcap V(D_{0k})|=1$ and $|S\bigcap V(D_{0l})|=1$ for distinct $i,j, k, l\in\{1,2,\cdots, 2^{n-1}\}$.

Without loss of generality, let $|S\bigcap V(D_{01})|=1, |S\bigcap V(D_{02})|=1, |S\bigcap V(D_{03})|=1$ and $|S\bigcap V(D_{04})|=1$. Let $x\in V(D_{01}), y\in V(D_{02}), z\in V(D_{03})$ and $w\in V(D_{04})$. See Figure $9$. Choose $n-1$ vertices $x_{1}, x_{2},\cdots, x_{n-1}$ from $D_{01}\setminus\{x\}$. Without loss of generality, let $x_{i}^{1}\in V(D_{1i})$, where $x_{i}^{1}$ is the outside neighbour of $x_{i}$ in $D_{n}$. Then choose $n-1$ vertices $y_{1}, y_{2},\cdots, y_{n-1}$ from $D_{02}$, $n-1$ vertices $z_{1}, z_{2},\cdots, z_{n-1}$ from $D_{03}$ and $n-1$ vertices $w_{1}, w_{2},\cdots, w_{n-1}$ from $D_{04}$ such that $y_{i}^{1}, z_{i}^{1}$ and $w_{i}^{1}$ all belong to $D_{1i}$, where $y_{i}^{1}, z_{i}^{1}$ and $w_{i}^{1}$ are the outside neighbour $y_{i}, z_{i}$ and $w_{i}$, respectively. Let $X=\{x_{1}, x_{2},\cdots, x_{n-1}\}, Y=\{y_{1}, y_{2},\cdots, y_{n-1}\}, Z=\{z_{1}, z_{2},\cdots, z_{n-1}\}$ and $W=\{w_{1}, w_{2},\cdots, w_{n-1}\}$. By Lemma~\ref{lem6}, there are $n-1$ internally disjoint paths $P_{1}, P_{2}, \cdots, P_{n-1} $ from $x$ to $X$ in $D_{01}$, $n-1$ internally disjoint paths $Q_{1}, Q_{2}, \cdots, Q_{n-1} $ from $y$ to $Y$ in $D_{02}$, $n-1$ internally disjoint paths $R_{1}, R_{2}, \cdots, R_{n-1} $ from $z$ to $Z$ in $D_{03}$ and $n-1$ internally disjoint paths $S_{1}, S_{2}, \cdots, S_{n-1} $ from $w$ to $W$ in $D_{04}$, respectively. It is possible that one of the paths $Q_{i}s(R_{i}s, S_{i}s)$ is a single vertex. As $D_{1i}$ is connected, there is tree $\widehat{T}_{i}$ connecting $x_{i}^{1}, y_{i}^{1}$, $z_{i}^{1}$ and $w_{i}^{1}$ in $D_{1i}$ for each $i\in\{1,2,\cdots, n-1\}$. Let $T_{i}=P_{i}\bigcup Q_{i}\bigcup R_{i}\bigcup S_{i}\bigcup \widehat{T}_{i} \bigcup x_{i}x_{i}^{1}\bigcup y_{i}y_{i}^{1}\bigcup z_{i}z_{i}^{1}\bigcup w_{i}w_{i}^{1}$ for each $i\in \{1, 2, \cdots, n-1\}$. Then $n-1$ internally disjoint $S$-trees $T_{i}$s for $1\leq i\leq n-1$ are obtained in $D_{n}$.
\begin{figure}[!ht]
\begin{center}
%\vskip4cm
\includegraphics[scale=0.7]{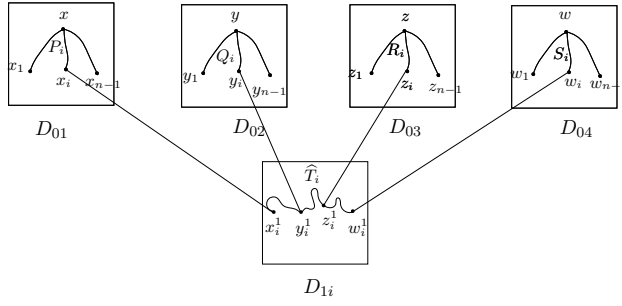}
\end{center}
\vskip0.5cm
\caption{The illustration of Case $1$}\label{F9}
\end{figure}

Case $2$. $|S\bigcap V(D_{0i})|=1, |S\bigcap V(D_{0j})|=1, |S\bigcap V(D_{0k})|=1$ and $|S\bigcap V(D_{1l})|=1$ for distinct $i,j, k$ and $i,j, k ,l \in\{1,2,\cdots, 2^{n-1}\}$.

Without loss of generality, let $|S\bigcap V(D_{01})|=1, |S\bigcap V(D_{02})|=1, |S\bigcap V(D_{03})|=1$ and $|S\bigcap V(D_{11})|=1$. Let $x\in V(D_{01}), y\in V(D_{02}), z\in V(D_{03})$ and $w\in V(D_{11})$. See Figure $10$. Choose $n-1$ vertices $x_{1}, x_{2},\cdots, x_{n-1}$ from $D_{01}$. Without loss of generality, let $x_{i}^{1}\in V(D_{1(i+3)})$, where $x_{i}^{1}$ is the outside neighbour of $x_{i}$. Then choose $n-1$ vertices $y_{1}, y_{2},\cdots, y_{n-1}$ from $D_{02}$ and $n-1$ vertices $z_{1}, z_{2},\cdots, z_{n-1}$ from $D_{03}$ such that $y_{i}^{1}, z_{i}^{1}\in V(D_{1(i+3)})$, where $y_{i}^{1}$ and $z_{i}^{1}$ are the outside neighbour $y_{i}$ and $z_{i}$ in $D_{n}$, respectively. Then choose $n-1$ vertices $w_{1}, w_{2},\cdots, w_{n-1}$ from $D_{11}$ such that $w_{i}^{0}$, the outside neighbour of $w_{i}$, does not belong to $D_{01}, D_{02}$ and $D_{03}$. Without loss od generality, let $w_{i}^{0}\in V(D_{0(i+3)})$. This can be done as $2^{n-1}> n+2$ for $n\geq 4$. Let $X=\{x_{1}, x_{2},\cdots, x_{n-1}\}, Y=\{y_{1}, y_{2},\cdots, y_{n-1}\}, Z=\{z_{1}, z_{2},\cdots, z_{n-1}\}$ and $W=\{w_{1}, w_{2},\cdots, w_{n-1}\}$. By Lemma~\ref{lem6}, there are $n-1$ internally disjoint paths $P_{1}, P_{2}, \cdots, P_{n-1} $ from $x$ to $X$ in $D_{01}$, $n-1$ internally disjoint paths $Q_{1}, Q_{2}, \cdots, Q_{n-1} $ from $y$ to $Y$ in $D_{02}$, $n-1$ internally disjoint paths $R_{1}, R_{2}, \cdots, R_{n-1} $ from $z$ to $Z$ in $D_{03}$ and $n-1$ internally disjoint paths $S_{1}, S_{2}, \cdots, S_{n-1} $ from $w$ to $W$ in $D_{11}$, respectively. It is possible that one of the paths $P_{i}s(Q_{i}s, R_{i}s, S_{i}s)$ is a single vertex. By Theorem~\ref{lem8}, $D_{n}[V(D_{1(i+3)})\bigcup V(D_{0(i+3)})]$ is connected. Then there is tree $\widehat{T}_{i}$ connecting $x_{i}^{1}, y_{i}^{1}$, $z_{i}^{1}$ and $w_{i}^{0}$ in $D_{n}[V(D_{1(i+3)})\bigcup V(D_{0(i+3)})]$. Let $T_{i}=P_{i}\bigcup Q_{i}\bigcup R_{i}\bigcup S_{i}\bigcup \widehat{T}_{i}\\\bigcup x_{i}x_{i}^{1}\bigcup y_{i}y_{i}^{1}\bigcup z_{i}z_{i}^{1}\bigcup w_{i}w_{i}^{0}$ for each $i\in \{1, 2, \cdots, n-1\}$. Then $n-1$ internally disjoint $S$-trees $T_{i}$s for $1\leq i\leq n-1$ are obtained in $D_{n}$.

\begin{figure}[!ht]
\begin{center}
%\vskip4cm
\includegraphics[scale=0.7]{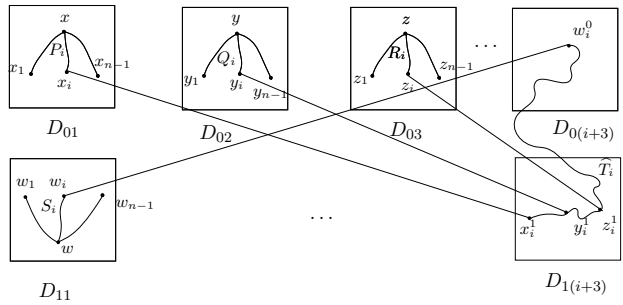}
\end{center}
\vskip0.5cm
\caption{The illustration of Case $2$}\label{F10}
\end{figure}

Case $3$. $|S\bigcap V(D_{0i})|=1, |S\bigcap V(D_{0j})|=1, |S\bigcap V(D_{1k})|=1$ and $|S\bigcap V(D_{1l})|=1$ for $i\neq j$, $k\neq l$ and $i,j, k, l\in\{1,2,\cdots, 2^{n-1}\}$.

Without loss of generality, let $|S\bigcap V(D_{01})|=1, |S\bigcap V(D_{02})|=1, |S\bigcap V(D_{11})|=1$ and $|S\bigcap V(D_{12})|=1$. Let $x\in V(D_{01}), y\in V(D_{02}), z\in V(D_{11})$ and $w\in V(D_{12})$. See Figure $11$. Choose $n-1$ vertices $x_{1}, x_{2},\cdots, x_{n-1}$ from $D_{01}$. Without loss of generality, let $x_{i}^{1}\in V(D_{1(i+2)})$, where $x_{i}^{1}$ is the outside neighbour of $x_{i}$ for $1\leq i\leq n-1$. Then choose $n-1$ vertices $y_{1}, y_{2},\cdots, y_{n-1}$ from $D_{02}$ such that $y_{i}^{1}\in V(D_{1(i+2)})$, where $y_{i}^{1}$ is the outside neighbour of $y_{i}$. Similarly, choose $n-1$ vertices $z_{1}, z_{2},\cdots, z_{n-1}$ from $D_{11}$. Without loss of generality, let $z_{i}^{0}\in V(D_{0(i+2)})$, where $z_{i}^{0}$ is the outside neighbour of $z_{i}$ for $1\leq i\leq n-1$. Then choose $n-1$ vertices $w_{1}, w_{2},\cdots, w_{n-1}$ from $D_{12}$ such that $w_{i}^{0}\in V(D_{0(i+2)})$, where $w_{i}^{0}$ is the outside neighbour of $w_{i}$. Let $X=\{x_{1}, x_{2},\cdots, x_{n-1}\}, Y=\{y_{1}, y_{2},\cdots, y_{n-1}\}, Z=\{z_{1}, z_{2},\cdots, z_{n-1}\}$ and $W=\{w_{1}, w_{2},\cdots, w_{n-1}\}$. By Lemma~\ref{lem6}, there are $n-1$ internally disjoint paths $P_{1}, P_{2}, \cdots, P_{n-1} $ from $x$ to $X$ in $D_{01}$, $n-1$ internally disjoint paths $Q_{1}, Q_{2}, \cdots, Q_{n-1} $ from $y$ to $Y$ in $D_{02}$, $n-1$ internally disjoint paths $R_{1}, R_{2}, \cdots, R_{n-1} $ from $z$ to $Z$ in $D_{11}$ and $n-1$ internally disjoint paths $S_{1}, S_{2}, \cdots, S_{n-1} $ from $w$ to $W$ in $D_{12}$, respectively. It is possible that one of the paths $P_{i}s(Q_{i}s, R_{i}s, S_{i}s)$ is a single vertex. By Theorem~\ref{lem8}, $D_{n}[V(D_{0(i+2)})\bigcup V(D_{1(i+2)})]$ is connected for each $i\in\{1,2,\cdots, n-1\}$. Then there is tree $\widehat{T}_{i}$ connecting $x_{i}^{1}, y_{i}^{1}$, $z_{i}^{0}$ and $w_{i}^{0}$ in $D_{n}[V(D_{0(i+2)})\bigcup V(D_{1(i+2)})]$. Let $T_{i}=P_{i}\bigcup Q_{i}\bigcup R_{i}\bigcup S_{i}\bigcup \widehat{T}_{i} \bigcup x_{i}x_{i}^{1}\bigcup y_{i}y_{i}^{1}\bigcup z_{i}z_{i}^{0}\bigcup w_{i}w_{i}^{0}$ for each $i\in \{1, 2, \cdots, n-1\}$. Then $n-1$ internally disjoint $S$-trees $T_{i}$s for $1\leq i\leq n-1$ are obtained in $D_{n}$.

\begin{figure}[!ht]
\begin{center}
%\vskip4cm
\includegraphics[scale=0.7]{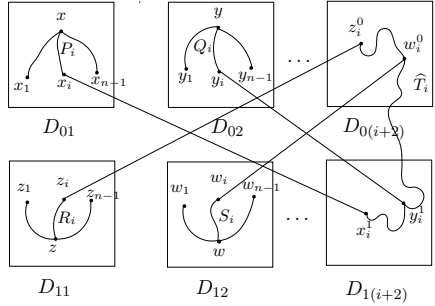}
\end{center}
\vskip0.5cm
\caption{The illustration of Case $3$}\label{F3}
\end{figure}

\hfill\qed

\begin{theorem}\label{thm1}
Let $D_{n}$ be an $n$-dimensional dual cube, then $\kappa_{4}(D_{n})=n-1$ for $n\geq 4$.
\end{theorem}

\f {\bf Proof.} As $D_{n}$ is $n$ regular, by Lemma~\ref{lem5}, $\kappa_{4}(D_{n})\leq \delta-1= n-1.$ To prove the result, we just need to show that $\kappa_{4}(D_{n})\geq n-1.$

Let $D_{01}, D_{02},\cdots, D_{02^{n-1}}$ be the clusters of class $0$ and $D_{11}, D_{12},\cdots$ and $D_{12^{n-1}}$ be the clusters of class $1$ of $D_{n}$. Let $x, y, z$ and $w$ be any four distinct vertices of $D_{n}$. For convenience, let $S=\{x, y, z, w\}.$ By the symmetry of $D_{n}$, we prove the result by considering the following cases.

Case 1. $x, y, z$ and $w$ belong the same cluster of $D_{n}$.

Without loss of generality, let $S\subseteq V(D_{01})$. Recall that $D_{01}$ is a copy of $Q_{n-1}$. By Theorem~\ref{lem3}, $\kappa_{4}(Q_{n-1})=n-2$. Then there are $(n-2)$-internally disjoint $S$-trees $T_{1}, T_{2},\cdots, T_{n-2}$ in $D_{01}$. Let $x^{1}, y^{1}, z^{1}$ and $w^{1}$ be the outside neighbours of $x, y, z$ and $w$ in $D_{n}$. By Lemma~\ref{lem1}, $\{x^{1}, y^{1}, z^{1}, w^{1}\}\subseteq V(D_{n}\setminus D_{01})$. By Theorem~\ref{lem8}, $D_{n}\setminus D_{01}$ is connected. Thus, there is a tree $\widehat{T}_{n-1}$ connecting $x^{1}, y^{1}, z^{1}$ and $w^{1}$ in $D_{n}\setminus D_{01}$. Let $T_{n-1}=\widehat{T}_{n-1}\bigcup xx^{1}\bigcup yy^{1}\bigcup zz^{1}\bigcup ww^{1}$, then $T_{1}, T_{2},\cdots, T_{n-1}$ are $(n-1)$-internally disjoint $S$-trees in $D_{n}$ and the result is desired.

Case 2. $x, y, z$ and $w$ belong to two different clusters of $D_{n}$.

By Theorem~\ref{Lem81}, $(n-1)$-internally disjoint $S$-trees $T_{1}, T_{2},\cdots, T_{n-1}$ can be obtained in $D_{n}$.

Case 3. $x, y, z$ and $w$ belong to three different clusters of $D_{n}$.

By Theorem~\ref{Lem82}, $(n-1)$-internally disjoint $S$-trees $T_{1}, T_{2},\cdots, T_{n-1}$ can be obtained in $D_{n}$.

Case 4. $x, y, z$ and $w$ belong to four different clusters of $D_{n}$.

By Theorem~\ref{Lem83}, $(n-1)$-internally disjoint $S$-trees $T_{1}, T_{2},\cdots, T_{n-1}$ can be obtained in $D_{n}$.

Thus, $\kappa_{4}(D_{n})\geq n-1$ and the result is desired.
\hfill\qed

By Theorem~\ref{thm1} and Lemma~\ref{lem4}, we have the following corollary.

\begin{cor}\label{cor1}
Let $D_{n}$ be an $n$-dimensional dual cube. Then $\kappa_{3}(D_{n})=n-1$ for $n\geq 4$.
\end{cor}

\section{The $r$-component connectivity of dual cubes}

In this section, we prove the $r$-component connectivity $c\kappa_{r}(D_{n})$ of $D_{n}$. To prove the result, the following result is useful.

\begin{lem}{\rm(\cite{Zh})}\label{lem9}
Any two vertices in $V(Q_{n})$ have exactly two common neighbours for $n\geq 3$ if they have any.
\end{lem}

\begin{theorem}{\rm(\cite{Ar})}\label{lem10}
Let $n$ and $k$ be positive integers such that $n\geq 3$ and $1\leq k \leq n-1$. Let $D_{n}$ be a dual cube network of order $n$ and $T\subseteq V(D_{n})$. If $|T|\leq kn-\frac{k(k+1)}{2}$, then $D_{n}-T$ is either connected or it has a large connected component and small components with at most $k-1$ vertices in total. Moreover, there is a set of vertices $T$ in $D_{n}$ such that $|T|= kn-\frac{k(k+1)}{2}+1$ and $D_{n}-T$ has a component containing exactly $k$ vertices.
\end{theorem}

Following, we show the $r$-component connectivity of the dual cube $D_{n}$.

\begin{theorem}\label{thm2}
Let $D_{n}$ be the $n$-dimensional dual cube, then $c\kappa_{r+1}(D_{n})=rn-\frac{r(r+1)}{2}+1$ for $n\geq 2$ and $1\leq r \leq n-1$.
\end{theorem}

\f {\bf Proof.} By Theorem~\ref{lem10}, $c\kappa_{r+1}(D_{n})\geq rn-\frac{r(r+1)}{2}+1$. Following, we show that $c\kappa_{r+1}(D_{n})\leq rn-\frac{r(r+1)}{2}+1$. Let $u\in V(D_{01})$, then $u$ has $n-1$ neighbours $u_{1}, u_{2},\cdots,$ and $u_{n-1}$ in $D_{01}$ as $D_{01}$ is a copy of $Q_{n-1}$. As $D_{n}$ is bipartite, then $u_{i}$ and $u_{j}$ are non-adjacent for $i\neq j$ and $i,j\in\{1,2,\cdots,n-1\}$. Let $1\leq r\leq n-1$ and $S=N_{D_{n}}(\{u_{1},u_{2},\cdots,u_{r}\})$. Then $u\in S$. By Lemma~\ref{lem9}, $|S|=r(n-1)-\frac{r(r+1)}{2}+1+r=rn-\frac{r(r+1)}{2}+1$. It is obvious that $D_{n}-S$ contains at least $r+1$ components, where at least $r$ of them are singletons. Hence, $c\kappa_{r+1}(D_{n})\leq rn-\frac{r(r+1)}{2}+1$. Combining with $c\kappa_{r+1}(D_{n})\geq rn-\frac{r(r+1)}{2}+1$, we have that $c\kappa_{r+1}(D_{n})=rn-\frac{r(r+1)}{2}+1$ for $n\geq 2$ and $1\leq r \leq n-1$.
\hfill\qed

\section{Concluding remarks}
The dual cube has some attractive properties to design interconnection networks. In this paper, we focus on the dual cube, which is an invariant of the hypercube and denoted by $D_{n}$. We first show that $\kappa_{4}(D_{n})=n-1$ for $n\geq 4$. That is, for any four vertices in $D_{n}$, there are $n-1$ internally disjoint trees connecting them in $D_{n}$.
As a corollary, we obtain $\kappa_{3}(D_{n})=n-1$ for $n\geq 4$. Furthermore, we show that $c\kappa_{r+1}(D_{n})=rn-\frac{r(r+1)}{2}+1$ for $n\geq 2$ and $1\leq r \leq n-1$. So far, the results about the generalized $k$-connectivity of networks are almost about $k=3$ and there are few results about larger $k$. We are interested in this topic and we would like to study the generalized $k$-connectivity of the dual cube for $k\geq 5$ in the future. In addition, we conjecture that $\kappa_{r}(D_{n})=n-1$ for $5\leq r \leq n $ and $n\geq 2$. We would like to prove the result for the dual cube and show the corresponding results for other classical networks such as arrangement graphs, balanced graphs, $(n, k)$-star networks and so on.

\section*{Acknowledgments}
This work was supported by the National Natural Science Foundation of China (Nos.11371052, 11731002, 11571035), the Fundamental Research Funds for the Central Universities (No. 2016JBM071, 2016JBZ012) and the $111$ Project of China (B16002).

\end{document}